\documentclass[runningheads]{lncse}
\usepackage{multicol} % needed for the author index
\usepackage{makeidx}  % allows for index generation
\def\Z{\bbbz}           % LNCSE code for whole numbers

\def\C{\bbbc}           % LNCSE code for complex numbers

\def\R{\bbbr}           % LNCSE code for real numbers

\newcommand{\D}{\partial}
\newcommand{\F}{\frac}
\newcommand{\Lp}{\left(}
\newcommand{\Rp}{\right)}
\newcommand{\Sum}[2]{\sum\limits_{#1=1}^{#2}}
\newcommand{\Sumo}[2]{\sum\limits_{#1}^{#2}}

\makeindex
\begin{document}
\mainmatter              % start of the contribution
\title{Cohomology of Lie Superalgebras of
       Hamiltonian Vector Fields: Computer Analysis}
%\subtitle{Computer Analysis}
%
\titlerunning{Cohomology of Lie Superalgebras}
% abbreviated title (for running head)
%
\author{Vladimir V. Kornyak}
\authorrunning{Vladimir V. Kornyak}  % abbreviated author list (for running head)
%
%%%% modified list of authors for the TOC (add the affiliations)
%\tocauthor{Vladimir V. Kornyak {\rm(Joint Institute for Nuclear Research)}}
%
\institute{
       Laboratory of Computing Techniques and Automation\\
       Joint Institute for Nuclear Research\\
       141980 Dubna, Russia}

\maketitle
\index{Kornyak@Vladimir V. Kornyak}
\begin{abstract}
In this paper we present the results of computation of
cohomology for some Lie (super)algebras of Hamiltonian
vector fields and related algebras. At present, the full
cohomology rings for these algebras are not known even
for the low dimensional vector fields. The partial
``experimental" results may give some hints for solution
of the whole problem.
The computations have been carried out
with the help of recently written program in C language.
Some of the presented results are new.
\end{abstract}
\section{Introduction and Basic Definitions}
There are many applications of the Lie (super)algebra
cohomology in mathematics: characteristic classes of foliations;
invariant differential operators; MacDonald-type combinatorial
identities, etc. (see~\cite{Fuks} for details).
Moreover, the cohomology is widely used in mathematical and theoretical
physics: construction of the central ex\-ten\-si\-ons and
deformations for Lie super\-al\-gebras;
construction of supergravity equations for $N$-extended Minkowski
super\-spa\-ces and search for possible models for these super\-spa\-ces;
study of stability for nonholonomic systems like ballbearings, gyroscopes,
electro-mechanical devices, waves in plasma, etc.; description of an analogue
of the curvature tensor for nonlinear nonholonomic constraints~\cite{GL};
new methods for the study of integrabi\-li\-ty of dynamical systems;
construction of the so-called {\em higher order Lie algebras}~\cite{Fil} which
allow in turn to construct the {\em Nambu mechanics} generalizing the ordinary
Hamiltonian mechanics~\cite{Takh};
construction of possible invariant effective actions of Wess-Zumino-Witten
type and the study of anomalies~\cite{D'Hoker}.
In~\cite{KonMan} the {\em cohomological field theories} have been aaplied
to enumerative problems of algebraic geometry.

General definitions and properties of cohomology of Lie algebras
and superalgebras are described in~\cite{Fuks}.
Let us recall briefly some basic definitions.

A {\em Lie superalgebra\/} is a ${\Z}_2$-graded algebra over a commutative
ring $K$ with a unit:
$$L = L_{\bar{0}} \oplus L_{\bar{1}}, \
u \in L_\alpha, \  v \in L_\beta, \ \alpha, \beta \in {\Z}_2 =
\{\bar{0}, \bar{1}\}\Longrightarrow [u,v] \in L_{\alpha + \beta}$$
The elements of $L_{\bar{0}}$ and $L_{\bar{1}}$ are called {\em even\/}
and {\em odd}, respectively.
We shall assume $K$ is one of the fields $\C$
or  $\R$. By definition, the {\em Lie product\/} (shortly, {\em bracket\/})
 $[\cdot ,\cdot ]$ satisfies the following axioms
\begin{eqnarray*}
& [u,v] = -(-1)^{p(u)p(v)} [v,u], &  {\qquad \mbox{skew-symmetry},} \label{ss_s} \\
& [u,[v,w]] = [[u,v],w] +(-1)^{p(u)p(v)} [v,[u,w]], &
{\qquad \mbox{Jacobi\ identity},}
\end{eqnarray*}
where $p(a)$ is the parity of element $a \in L_{p(a)}$.

A {\em module\/} over a Lie superalgebra $A$ is a vector space $M$
(over the same field $K$) with a mapping $A \times M \to M$, such that
$[a_1, a_2] m = a_1 (a_2 m) - (-1)^{p(a_1)p(a_2)} a_2 (a_1 m)$,
where $a_1, a_2 \in A$, $m \in M$.
The most important for our purposes are {\em trivial\/} ($M$ is an arbitrary
vector space, e.g.,$M = K$; $am = 0$), {\em adjoint\/} ($M = A;
am = [a,m]$) and {\em coadjoint\/}  ($M = A'; am = \{a,m\}$ is coadjoint action)
modules.

A {\em cochain complex\/} is a sequence of linear spaces $C^k$
with linear mappings $d^k$
\begin{equation}
0\to C^0\stackrel{d^0}{\longrightarrow}\cdots
\stackrel{d^{k-2}}{\longrightarrow}C^{k-1}\stackrel{d^{k-1}}
{\longrightarrow} C^k\stackrel{d^k}{\longrightarrow}C^{k+1}
\stackrel{d^{k+1}}{\longrightarrow}\cdots,
\label{cocomplex}
\end{equation}
where the linear space $C^k = C^k(A;M)$ is a super skew-symmetric
$k$-linear mapping $A \times \cdots \times A \to M$, $C^0 = M$ by definition.
The super skew-symmetry means symmetry w.r.t. transpositions of odd adjacent
elements of $A$ and antisymmetry for all other transpositions of adjacent elements.
Elements of $C^k$ are called {\em cochains}.

The linear mapping $d^k$ (or, briefly, $d$) is called the {\em differential\/}
and satisfies the following property: $d^k \circ d^{k-1} = 0$ (or $d^2 = 0$).

The cochains mapped into zero by the differential are called the
{\em cocycles}, i.e., the space of cocycles is
$$Z^k = {\rm Ker}\ d^k = \{C^k \; | \; dC^k = 0\}.$$

The cochains which can be represented as differentials
of other cochains are called the {\em coboundaries}, i.e., the space
of coboundaries is
$$B^k = {\rm Im}\ d^{k-1} = \{C^k \; | \; C^k = dC^{k-1}\}.$$
Any coboundary is obviously a cocycle.

The non-trivial cocycles, i.e., those which are not coboundaries,
form the {\em cohomology}. In other words, the cohomology is
the quotient space $$H^k(A;M) = Z^k/B^k.$$

The explicit form of the differential for a Lie superalgebra is
\begin{eqnarray*}
& dC(e_0,\ldots,e_q;O_{q+1},\ldots,O_k)\ = & \\
&  \Sumo{i<j}{q}(-1)^j C(e_0,\ldots,e_{i-1},[e_i,e_j],\ldots,
    \widehat{e_j},\ldots,e_q;O_{q+1},\ldots,O_k)\ + &\\
& (-1)^{q+1}\Sumo{i=0}{q}\Sumo{j=q+1}{k}C(e_0,\ldots,e_{i-1},[e_i,O_j],\ldots,
   e_q;O_{q+1},\ldots,\widehat{O_j},\ldots, O_k)\ + &\\
& (-1)^{i+1}\Sumo{i=q+1}{k-1}\Sumo{j=q+2}{k}C(e_0,\ldots,e_q;O_{q+1},\ldots,
   O_{i-1},[O_i,O_j],\ldots,\widehat{O_j},\ldots, O_k) &\\
& + \Sumo{i=0}{q}(-1)^{i+1} e_i
   C(e_0,\ldots,\widehat{e_i},\ldots,e_q;O_{q+1},\ldots,O_k) &\\
& +\ (-1)^q \Sumo{i=q+1}{k} O_i
   C(e_0,\ldots,e_q;O_{q+1},\ldots,\widehat{O_i},\ldots,O_k). &
\end{eqnarray*}
Here $e_i$ and $O_i$ are even and odd elements of the algebra, respectively,
and the hat ``$\;\widehat{\ }\;$" marks the omitted elements.

Here are some properties and statements we use in the sequel.

An algebra and a module are called {\em graded\/} if they can be
presented as sums of homogeneous components in a way compatible
with the algebra bracket and the action of the algebra on the module:
$$A = \oplus_{g \in G}\, A_g, \; M = \oplus_{g \in G}\, M_g, \;
[A_{g_1}, A_{g_2}] \subset A_{g_1+g_2}, \;
A_{g_1}M_{g_2} \subset M_{g_1+g_2},$$
where $G$ is some abelian (semi)group. We assume $G = \Z$ in this paper.
To avoid confusion, we use in the sequel the terms {\em grade} and {\em degree}
for element of $G$ and number of cochain arguments, respectively.
The grading in the algebra and module induces a grading on cochains
and, hence, in the cohomology:
$$C^{*}(A;M) = \oplus_{g \in G}\, C_g^{*}(A;M), \quad
 H^{*}(A;M) = \oplus_{g \in G}\, H_g^{*}(A;M).$$
This property allows one to compute the cohomology separately
for different ho\-mo\-ge\-ne\-ous components; this is especially useful when
the homogeneous components are finite--dimensional.

If there is an element $a_0 \in A$, such that eigenvectors of the operator
$a \mapsto [a_0,a]$ form a (topological) basis of algebra $A$,
then $H^{*}(A) \simeq H_0^{*}(A).$ In other words, all the non-trivial
cocycles of the cohomology in the trivial module lie in the zero grade component.
The element $a_0$ is called an {\em internal grading element}.
If also eigenvectors of the operator $m \mapsto a_0m$
form a topological basis of module $M$, then the same statement holds
for the cohomology in the module $M$: $H^{*}(A;M) \simeq H_0^{*}(A;M).$

In the case of trivial module, the exterior multiplication of cochains
provides the cohomology with a structure of graded ring. There are
also another multiplicative structures in cohomology, but we shall not
use them in this work.

\section{Outline of Algorithm and Its Implementation}
\noindent
To compute the cohomology one needs to solve the equation
\begin{equation}
dC^k = 0,
\label{dck}
\end{equation}
and throw away those solutions of (\ref{dck}) which can be
expressed in the form
$$C^k = dC^{k-1}.$$
In some exceptional cases it is possible to solve equation (\ref{dck})
in closed form.
Often, in the case of Lie superalgebras of vector fields,
determining equation (\ref{dck}) is a system of linear homogeneous
functional equations with integer arguments. Unfortunately
there is no general method for solving
such systems in closed form. Hence, we need to carry out the corresponding
computation ``numerically". There are several packages for computing
cohomology of Lie algebras and superalgebras written in
{\em Reduce\/}~\cite{LP},~\cite{PH} and {\em Mathematica\/}~\cite{GL}.
Some new results
were obtained completely or partially with the help of these packages.
However, these packages, being based on general purpose computer algebra
systems, appeared to be too inefficient for large real problems.
In view of this, we wrote the program in C language~\cite{Korn98}.

The C code, of total length near 14000 lines, contains about 300
functions
realizing top level algorithms, simplification of indexed objects,
working
with Grassmannian objects, exterior calculus, linear algebra,
substitutions, list processing, input and output, etc.
As internal structures we use 8 types of lists for different objects.
We represent Grassmann monomials by integer numbers using one-to-one
correspondence between (binary codes of) non-negative integers and
Grassmann monomials. This representation allows one efficiently to
implement
the operations with Grassmann monomials by means of the basic computer
commands.

The program performs sequentially the following steps:
\begin{enumerate}
\item {\em Reading input information.}
\item {\em Constructing a basis\/} for the algebra. The basis can be
 read from the input file; otherwise the program constructs it from
 the definition of the algebra.
 Non-trivial computations
 at this step arise only in the case of divergence-free algebras.
 The basis elements of such algebras should satisfy some conditions.
 In fact, we should construct the basis elements of a subspace given by
 a system of linear equations. The task is thereby reduced to some
 problem of linear
 algebra combined with shifts of indices. For example, among the
 divergence-free
 conditions for the {\em special Buttin} algebra ${\mathrm{SB}(3)}$
 there are the
 following two equations
 $$ia_{ijk; UV} -(k+1)a_{i-1,j,k+1; VW} = 0,$$
 $$ia_{ijk; UW} +(j+1)a_{i-1,j+1,k; VW} = 0.$$
 Here $a_{ijk; UV},\ldots$ are coefficients at the monomials
 $p^iq^jr^kUV,\ldots$
 in the generating function; $p, q, r$ and $U, V, W$ are even and odd
 variables,
 respectively. First of all, we have to shift indices $j$ and $k$ in
 the second
 equation to reduce the last terms of both equations to the same
 multiindices.
 Then, using some simple tricks of linear algebra, we can easily construct the
 corresponding basis element
 $$E_{ijk} = (k+1)p^iq^jr^kUV - jp^iq^{j-1}r^{k+1}UW + ip^{i+1}q^jr^{k+1}VW.$$

\item {\em Constructing the commutator table\/} for the algebra
 (if this table has not been read from the input file).
\item {\em Creating the general form\/} of expressions for coboundaries
 and determining equa\-ti\-ons for cocycles.
\item {\em Transition to a particular grade\/} in general expressions.
 At this step expressions for coboundaries take the form
 ${\bf x = bt}$, equations
 for cocycles take the form ${\bf Zx = 0}$, where vector ${\bf x}$
 corresponds to $C^k$, parameter vector ${\bf t}$ corresponds
 to $C^{k-1}$, matrices ${\bf Z, b}$
 correspond to the differential $d$. All these vector spaces are
 finite-dimensional for any particular grade.
\item {\em Computing the quotient space $H^k(A;M) = Z^k/B^k$}.
 Here cocycle subspace $Z^k$
 is given by relations ${\bf Zx = 0}$, and coboundary subspace $B^k$ is given
 parametrically by ${\bf x = bt}$.

 Substeps:
 \begin{enumerate}
 \item Eliminate ${\bf t}$ from  ${\bf x = bt}$ to get equations
  ${\bf Bx = 0}$
 \item Reduce both relations ${\bf Bx = 0}$ and ${\bf Zx = 0}$ to
  the canonical
  form by Gauss elimination. If ${\rm rank} {\bf B} = {\rm rank} {\bf Z}$,
  then there is no non-trivial cocycle;
  otherwise go to Substep (c).
 \item Set ${\bf Bx = y}$ and substitute these relations into ${\bf Zx = 0}$
  to get relations ${\bf Ay = 0}$. The {\em parametric} (non-leading)
  $y'$s of the last relations are non-trivial cocycles; that is, they
  form a basis of the cohomology.
 \end{enumerate}
 In fact, the above procedure is based on the relation for quotient
 spaces
 $$Z/B = \frac{Y/B}{Y/Z},$$
 where $Y$ is an artificially introduced space combining the above $x'$s and $y'$s.
 \item {\em Output the non-trivial cocycles.}
\end{enumerate}

\section{Hamiltonian Vector Fields and Related Algebras}
\noindent
To define the formal vector fields on the supermanifold
of the superdimension $(2n|m)$ we consider the sets of even
$p_1,\ldots,p_n, q_1,\ldots,q_n;$
and odd (called also {\em Grassmann\/}) $U_1,\ldots,U_m$ variables,
and formal power series $f, g,\ldots,$ in these variables.
These power series are called {\em generating functions,}
because the vector fields considered in this work can be expressed
in terms of the derivatives of $f, g,\ldots$
The (super)commutator of vector fields induces the bracket on generating
functions.
The Lie superalgebra of {\em Poisson} vector fields
${\mathrm{Po}(2n|m)}$ is
a set of generating functions with the bracket
  $$\{f,g\}=\Sum{i}{n}\Lp\F{\D f}{\D p_i}\F{\D g}{\D q_i}
                            -\F{\D f}{\D q_i}\F{\D g}{\D p_i}\Rp
        -(-1)^{p(f)}\Sum{k}{m}\F{\D f}{\D U_k}\F{\D g}{\D U_k}.$$
The {\em Hamiltonian} superalgebra is a quotient algebra of the
Poisson algebra with respect to its center:
$${\mathrm{H}(2n|m) = \mathrm{Po}(2n|m)/Z}.$$
Observe that ${\mathrm{H}(0|m)}$ is not simple;
it has a simple ideal of
codimension 1 denoted by ${\mathrm{SH}(0|m)}$
and called {\em special Hamiltonian superalgebra}.
The algebra ${\mathrm{SH}(0|m)}$ contains the subalgebra
${\mathrm{O}(m)}$ and this fact can be used for analysis of
the structure of
the cohomology ring ${ H^*(\mathrm{SH}(0|m))}$.
Note that all the algebras depending only on the odd variables are
finite-dimensional. It does not mean however that their cohomologies
are finite-dimensional too.
All the above algebras are graded due to prescribed grading of
the variables.
The {\em standard\/} grading assumes all variables $q_i, p_i, U_i$ have
the grade 1.

\section{Computations}
\noindent
In the below tables and formulas we use the small $a,b,c,\ldots$
and capital $A,B,C,\ldots$ letters for even and odd cocycles,
respectively.
The optional superscript and subscript indicate the cochain degree and
grade, correspondingly. The letters without indices denote the genuine
cocycles, i.e., generating elements of the cohomology ring.
Empty position in the tables means the absence of non-trivial cocycles
in the given degree and grade.
The columns containing only trivial cocycles are omitted.
We use the notations $p,q$ for even and $U_i$ for odd
variables of the vector field generating functions.

\subsection{Special Hamiltonian Superalgebra}
\begin{table}[h!]
\caption{$H^n_g(\mathrm{SH}(0|4))$}
\begin{center}
\begin{tabular}{|l|c|c|c|c|c|c|c|}
\hline
$n \backslash g$ & -6 & -4 & -2 & 0 & 2 & 4 & 6 \\
\hline
1 & & & & & & & \\
\hline
2 & & & $a$ & $b$ & $c$ & & \\
\hline
3 & & & & $f$ & & & \\
\hline
4 & & $a^2$ & $ab$ & $b^2$ & $bc$ & $c^2$ & \\
\hline
5 & & & $af$ & $bf$ & $cf$ & & \\
\hline
6 & $a^3$ & $a^2b$ & $ab^2$ & $b^3$ & $ac^2$ & $bc^2$ & $c^3$ \\
\hline
\end{tabular}
\end{center}
\end{table}

In Table 1 the cohomology ring $H^*(\mathrm{SH}(0|4))$ is presented.
One can see that this ring is generated by four generators $a, b, c, f$
obeying to the relations $ac-b^2=0$ and $f^2=0$.
The explicit form of the generators given by the computer is
\begin{eqnarray*}
a&=&C(U_4,U_4)=C(U_1,U_1)=C(U_2,U_2)=C(U_3,U_3), \\
b&=&C(U_4,U_1 U_2 U_3)=C(U_1,U_2 U_3 U_4)=C(U_2,U_1 U_3 U_4)=C(U_3,U_1 U_2 U_4), \\
c&=&C(U_2 U_3 U_4,U_2 U_3 U_4)=C(U_1 U_2 U_3,U_1 U_2 U_3)=C(U_1 U_2 U_4,U_1 U_2 U_4) \\
 & &=C(U_1 U_3 U_4,U_1 U_3 U_4), \\
f&=&C(U_1 U_4,U_2 U_4,U_3 U_4)+\frac{1}{2}C(U_1 U_4,U_1,U_1 U_2 U_3)+ \\
 & &\frac{1}{2}C(U_1 U_4,U_4,U_2 U_3 U_4)+\frac{1}{2}C(U_2 U_4,U_2,U_1 U_2 U_3)- \\
 & &\frac{1}{2}C(U_2 U_4,U_4,U_1 U_3 U_4)+\frac{1}{2}C(U_3 U_4,U_3,U_1 U_2 U_3)+\\
 & &\frac{1}{2}C(U_3 U_4,U_4,U_1 U_2 U_4) = \ldots
\end{eqnarray*}
We have omitted for brevity the equivalent forms of generator $f$ in
the last formula.

\paragraph{Note 1.}
The cohomology of $\mathrm{SH}(0|4)$ has been computed for the first
time in~\cite{FL}
by D.Fuchs and D.Leites\hspace{-1mm}~\footnote{D.Leites informed us that
they missed the cocycle $f$ which was discovered
later by A. Shapovalov with the help of the program written
by P.Grozman.} by hand.
We present this example here as a rather short illustration demonstrating many
features of cohomology ring structure.
\paragraph{Note 2.}
An interesting approach based on algebraic geometry was sug\-ges\-ted by
C.~Gruson~\cite{G}. Her method enabled her to compute cohomology with
trivial coefficients of exceptional simple finite dimensional Lie
superalgebras. Though remarkably beautiful, it is not universal:
it fails
to work in various natural and interesting problems, e.g., for
$\mathrm{PSL}(n|n)=\mathrm{SL}(n|n)/Z$. It is unclear if Gruson's
method works for nontrivial coefficients. Observe that
$\mathrm{PSL}(2|2)=\mathrm{SH}(0|4)$,
whose cohomology we presented in Table 1.
\paragraph{}
\hspace{3mm}The structure of $H^*(\mathrm{SH}(0|m))$ has some peculiarities at $m = 4.$
Computations for $m = 3,5,6$ revealed two generators: even 2-cocycle
$$a=a^2_{-2}=C(U_1,U_1)=\ldots=C(U_m,U_m)$$
and 3-cocycle
$$\begin{array}{lllllll}
f & = & f^3_0 & = & C(U_1 U_2,U_1 U_3,U_2 U_3), & f^2=0, & m=3 \\
F & = & F^3_5 & = & C(U_1 U_3 U_4 U_5,U_2 U_3 U_4 U_5,U_3 U_4 U_5) = \ldots, & F^2=0, & m=5 \\
f & = & f^3_6 & = & C(U_3 U_4 U_5 U_6,U_4 U_5 U_6,U_1 U_2 U_4 U_5 U_6), & f^2=0, & m=6.
\end{array}$$
As we checked, there are no other generators in the case $m=3$
up to 16-cocycles.

It would be interesting to look how extensions of the algebra
influence on the structure of its cohomology ring. Tables 2 and 3
present the
cohomology structure for the superalgebras $\mathrm{H}(0|4)$ and
$\mathrm{Po}(0|4).$
Our consideration of the multiplicative structure for these cohomologies
is very preliminary.
In fact, there is a need to write a program for multiplication
and comparison of cocycles modulo coboundaries because corresponding computations
are rather tedious and error prone.
In Tables 2 and 3 the cocycle $a^5_{-2}$ is a linear combination of $af$ and $a^2b$,
$f'$ is cocycle and  $\alpha$ is a 2-cochain. The cocycle $e = C(U_1\ldots U_m)$
satisfies the relation $e^2=0$ for $m$ even  and is a free generator for $m$ odd.

\begin{table}[h!]
\caption{$H^n_g(\mathrm{H}(0|4))$}
\begin{center}
\begin{tabular}{|l|c|c|c|c|c|c|c|}
\hline
$n \backslash g$ & -6 & -4 & -2 & 0 & 2 & 4 & 6 \\
\hline
1 & & & & & $e$ & & \\
\hline
2 & & & $a$ & & & & \\
\hline
3 & & & & $f=f'+ae$ & & $r=e\alpha$ & \\
\hline
4 & & $a^2$ & & & $ef$ & & \\
\hline
5 & & & $a^5_{-2}$ & & & & $s=e\alpha^2$ \\
\hline
6 & $a^3$ & & & & & $rf$ & \\
\hline
\end{tabular}
\end{center}
\end{table}
\begin{table}[h!]
\caption{$H^n_g(\mathrm{Po}(0|4))$}
\begin{center}
\begin{tabular}{|l|c|c|c|c|}
\hline
$n \backslash g$ & 0 & 2 & 4 & 6 \\
\hline
1 & & $e$ & & \\
\hline
2 & $b$ & & & \\
\hline
3 & $f$ & & $r=e\alpha$ & \\
\hline
4 & & $ef, h = b\alpha$ & & \\
\hline
5 & $bf$ & & & $s=e\alpha^2$ \\
\hline
6 & & & $rf, k = b\alpha^2$ & \\
\hline
\end{tabular}
\end{center}
\end{table}

%\newpage
\subsection{Algebras $\mathrm{H}(2|0)$ and $\mathrm{Po}(2|0)$}

The case of supermanifold with even variables leads to infinite-dimensional
algebras and is much more difficult for the analysis. There are a few results
concerning the cohomology of Lie algebra $\mathrm{H}(2|0).$
In~\cite{GKF}
it has been proved that
 $\dim H^2(\mathrm{H}(2|0)) \geq 1, \dim H^5(\mathrm{H}(2|0)) \geq 1,
\dim H^7(\mathrm{H}(2|0)) \geq 2$ and $\dim H^{10}(\mathrm{H}(2|0))
\geq 1.$
In~\cite{Perchik} the inequality $\dim H^*(\mathrm{H}(2|0)) \geq 112$
has been obtained. The both works were based on the extraction of some easier
to handle
subcomplex of the full cohomological complex and application of the computer
analysis to this subcomplex. Besides, some facts about
$\dim H^q(\mathrm{H}(2k|0))$
for low degrees $q \leq k$ are known~\cite{GSH}.

Every 3--valent graph with oriented vertices, or every oriented rational
homology 3--sphere, can be associated with a cohomology class of
the Lie algebra
$\mathrm{H}(2k|0)$ for arbitrary $k$ (see~\cite{Ko1,BN,LT}).
This cohomology class called
{\em graph cohomology}~\footnote{This cohomology class can be calculated
via certain
complex constructed from finite graphs.} has been used for construction of
{\em Rozansky--Witten invariants}~\cite{Konts}.

Some cocycles from $H^*(\mathrm{Po}(2|0))$ (denoted by $p^n_g$) and
$H^*(\mathrm{H}(2|0))$
(denoted by $h^n_g$) obtained by the program are presented in Table 4.
We carried
out the computations up to degree 10 and grade 4.

\begin{table}[h!]
\caption{$H^n_g(\mathrm{Po}(2|0))$ and $H^n_g(\mathrm{H}(2|0))$}
\begin{center}
\begin{tabular}{|l|c|c|c|}
\hline
$n \backslash g$ & -4 & -2 & 0 \\
\hline
2 & & $h^2_{-2}$ & \\
\hline
3 & $p^3_{-4}$ & & \\
\hline
4 & & & \\
\hline
5 & & $p^5_{-2}, h^5_{-2}$ &  \\
\hline
6 & $p^6_{-4}$ & & \\
\hline
7 & & & $p^7_{0}, h^7_{0}$ \\
\hline
8 & & $p^8_{-2}$ & \\
\hline
\end{tabular}
\end{center}
\end{table}

%\newpage
\subsection{Algebras $\widehat{\mathrm{H}}(2|0)$ and
$\widehat{\mathrm{Po}}(2|0)$}

If we add the grading element $G$ to the algebra then the non-trivial cocycles lie
in zero grade only. In this case the space of cochains is finite-dimensional
and we can compute the full cohomology.
Thus, let's consider the algebras
$\widehat{\mathrm{Po}}(2|0) = \mathrm{Po}(2|0)\oplus \mathrm{Span}(G)$
and $\widehat{\mathrm{H}}(2|0) = \mathrm{H}(2|0)\oplus \mathrm{Span}(G).$
The cohomologies
$H^*(\widehat{\mathrm{Po}}(2|0))$
and $H^*(\widehat{\mathrm{H}}(2|0))$ are Grassmann algebras of
the superdimension $(2|2)$.
These algebras are generated by unit\footnote{Taking into account
the initial
part of cochain complex (\ref{cocomplex}) one can consider 1 formally as
a ``non-trivial" zero-cocycle.},  two
cocycles $a^1_0 = C(G)$ and $a^7_0$ and contain also 8-cocycle
$a^1_0a^7_0$.
For the case of Hamiltonian algebra the explicit form of $a^7_0$ is
$$a^7_0 = C(q,p,q^2,pq,p^2,q^3,p^3) - 3 C(q,p,q^2,pq,p^2,pq^2,p^2q).$$
In the Poisson case the expression for $a^7_0$ is much longer.

\section{Conclusion}
The computation of cohomology is a typical problem with
the combinatorial explosion. Nevertheless, some results
can be obtained with the help of computer having an efficient
enough program.
On the other hand, physicists are interested mainly in the
second cohomologies describing the central extensions and deformations.
Such cohomologies can be computed rather easily even for large algebras.
Some essential possibilities remain for increasing the efficiency of the program.
Besides, it would be useful also to write a separate program for investigating
the multiplicative structure of cohomology ring.

\section*{Acknowledgements}
I would like to thank D. Leites for initiating this work
and helpful communications.
I am also grateful to V. Gerdt and O. Khudaverdian for fruitful
discussions and useful advises.
This work was supported in part by INTAS project No. 96-184 and
RFBR project No. 98-01-00101.

%\newpage


\begin{thebibliography}{99}
\bibitem{LT}
Le Thang, An invariant of integral homology 3-spheres which is universal
for all finite type invariants, q--alg/9601002.
\end{thebibliography}
\end{document}